\documentclass[11pt,reqno]{amsart}
\usepackage{amssymb,amsmath}
\newcommand{\Z}{\mathbb{Z}}
\newcommand{\R}{\mathbb{R}}
\newcommand{\Q}{\mathbb{Q}}
\newcommand\Ha{\mathbb{H}}
\newcommand\ef{\mathcal{F}}
\renewcommand\theta{\vartheta}
\renewcommand\phi{\varphi}
\newcommand\sgn{\operatorname{sgn}}
\newcommand\smod[1]{\operatorname{mod} #1}
\newtheorem{theorem}{Theorem}[section]
\newtheorem{lemma}[theorem]{Lemma}
\newtheorem{remark}[theorem]{Remark}
\newtheorem{definition}[theorem]{Definition}
\numberwithin{equation}{section}

\title{Maass waveforms arising from $\sigma$ and related indefinite theta functions}
\date{\today}
\author{Sander Zwegers}
\address{School of Mathematical Sciences, University College Dublin, Belfield, Dublin 4, Ireland}
\email{sander.zwegers@ucd.ie}
\begin{document}
\begin{abstract}
In this paper we consider an example of a Maass waveform which was constructed by Cohen from a function $\sigma$, studied by Andrews, Dyson and Hickerson, and it's companion $\sigma^*$. We put this example in a more general framework.
\end{abstract}
\maketitle
\section{Introduction}
In \cite{ADH} Andrews, Dyson and Hickerson consider the following two $q$-hypergeometric functions, the first given in Ramanujan's ``Lost'' Notebook and the second, its ``partner'', discovered later
\begin{equation*}
\begin{split}
\sigma(q) &:= \sum_{n=0}^\infty \frac{q^{n(n+1)/2}}{(1+q) (1+q^2) \cdots (1+q^n)},\\
\sigma^*(q) &:= 2\sum_{n=1}^\infty \frac{(-1)^n q^{n^2}}{(1-q) (1-q^3) \cdots (1-q^{2n-1})}.
\end{split}
\end{equation*}
In that paper they find several identities for $\sigma$ and $\sigma^*$ which involve indefinite quadratic forms. Using these identities they are able to explain some amazing properties of the coefficients of these functions. Typical examples of identities they find are
\begin{equation}\label{sigma}
\begin{split}
q^{1/24} \sigma(q) &=\Bigl( \underset{n-j\geq 0}{\sum_{n+j\geq 0}} + \underset{n-j<0}{\sum_{n+j<0}}\Bigr) (-1)^{n+j}  q^{\frac{3}{2}(n+\frac{1}{6})^2-j^2},\\
q^{-1/24} \sigma^*(q) &=\Bigl( \underset{2j-3n> 0}{\sum_{2j+3n\geq 0}} + \underset{2j-3n\leq 0}{\sum_{2j+3n< 0}}\Bigr) (-1)^{n+j}  q^{-\frac{3}{2}(n+\frac{1}{6})^2+j^2}.
\end{split}
\end{equation}
These are equations (1.5) and (5.1) in \cite{ADH}, rewritten for the purpose of this paper. Note that the right hand sides can be viewed as theta functions attached to an indefinite quadratic form, or indefinite theta functions for short. Further note that the right hand side of the first equation closely resembles
\begin{equation*}
\Bigl( \underset{n-j\geq 0}{\sum_{n+j\geq 0}} - \underset{n-j<0}{\sum_{n+j<0}}\Bigr) (-1)^{n+j}  q^{\frac{3}{2}(n+\frac{1}{6})^2-j^2},
\end{equation*}
a function which is related to one of Ramanujan's sixth order mock theta functions. In \cite{zwegers} we find a general theory for indefinite theta functions of this second type, that is, sums of the form
\begin{equation*}
\sum_{\nu\in a+\Z^r} \left\{ \sgn(l_1(\nu))-\sgn(l_2(\nu))\right\} q^{Q(\nu)} e^{2\pi i l_3(\nu)},
\end{equation*}
where $Q$ is a quadratic form of signature $(r-1,1)$ on $\R^r$, $a\in\R^r$ and $l_1$, $l_2$ and $l_3$ are suitable linear functions. For special choices of $Q$, $a$ and $l_i$ this indefinite theta function is a holomorphic modular form of weight $r/2$ (on some subgroup of $\operatorname{SL}_2(\Z)$ with some multiplier system). In general, however, we get a more complicated object, namely the ``holomorphic part'' of a real-analytic modular form.

For $\sigma$ and $\sigma^*$ the relation with automorphic forms is more subtle: in \cite{cohen}, Cohen interprets the identities for $\sigma$ and $\sigma^*$ in terms of the theory of Maass waveforms. For this he defines the coefficients $T(n)$ by
\begin{equation*}
\underset{n\equiv 1\smod{24}}{\sum_{n\in\Z}} T(n) q^{|n|/24} = q^{1/24} \sigma (q) + q^{-1/24} \sigma^* (q),
\end{equation*}
and uses them to construct
\begin{equation*}
\phi_0(\tau) := y^{1/2} \underset{n\neq 0}{\sum_{n\in\Z}} T(n)\ e^{2\pi inx/24} K_0 (2\pi |n|y/24),
\end{equation*}
where $\tau=x+iy\in\Ha$ and $K_0$ is a modified Bessel function of the second kind. Cohen then shows that $\phi_0$ is a weight 0 Maass form: it transforms as a modular form of weight 0 on $\Gamma_0(2)$, with some (explicit) multiplier, but instead of being holomorphic on $\Ha$ the function satisfies
\begin{equation*}
\Delta_0 \phi_0 = \frac{1}{4} \phi_0,
\end{equation*}
where $\Delta_0 = -y^2 \left( \frac{\partial^2}{\partial x^2} + \frac{\partial^2}{\partial y^2}\right)=-4y^2 \frac{\partial^2}{\partial \tau \partial \overline{\tau}}$ is the weight 0 Laplace operator.

The differential equation satisfied by $\phi_0$ follows immediately from the differential equation satisfied by $K_0$, which is $\left(x \frac{\partial^2}{\partial x^2} + \frac{\partial}{\partial x} -x \right) K_0(x)=0$.
The hard part is to show that it also transforms like a modular function. For this he uses the Mellin transform and $L$-series.

The question that we address here is: what if we change the quadratic form and/or the linear restrictions in the right hand sides of equation \eqref{sigma} and do the same construction, will the resulting function still be a Maass waveform? The precise definition of the functions we consider is given in Definition \ref{defphi} below. From the construction it is immediately clear that these functions again are eigenfunctions of $\Delta_0$ with eigenvalue $1/4$, but do they also transform like modular functions? The answer is that in general they do not and we get interesting, more general, objects. In certain special cases we will again get a Maass waveform.

The outline of the paper is as follows: in the next section we will state some definitions and the main results. In section 3 we will give technical details for some of the definitions and in section 4 we will prove the main results. In section 5 we will consider the action of $\operatorname{Aut}^+(Q,\Z^2)$, which will be used in sections 6 and 7 to give some nice examples: in section 6 we consider Cohen's $\phi_0$ function and in section 7 we will obtain a family of cases where our construction gives a Maass waveform.

\section*{Acknowledgements}
The author wishes to thank Soon-Yi Kang and Ken Ono for suggesting this project.

\section{Definitions and statement of results}
Let $Q$ be a binary quadratic form of signature $(1,1)$. Let $A$ be the symmetric $2\times 2$-matrix such that $Q(\nu)= \frac{1}{2} \nu^t A \nu$. Throughout we assume that $A$ has integer coefficients. Further let $B$ be the associated bilinear form $B(\nu,\mu)=\nu^t A\mu=Q(\nu+\mu)-Q(\nu)-Q(\mu)$.

The set of vectors $c\in\R^2$ with $Q(c)=-1$ has two components. If $B(c_1,c_2)<0$ then $c_1$ and $c_2$ belong to the same component, while if $B(c_1,c_2)>0$ then $c_1$ and $c_2$ belong to opposite components. Let $C_Q$ denote one of the two components. If $c_0$ is in that component, then $C_Q$ is given by
\begin{equation*}
C_Q := \{ c\in \R^2 \mid Q(c)=-1\ \text{and}\ B(c,c_0)<0\}.
\end{equation*}
Since $Q$ is of signature $(1,1)$, it splits over $\R$ as the product of two linear functions, that is, we can write $Q(\nu)=Q_0(P\nu)$ with $P\in \operatorname{GL}_2(\R)$ and $Q_0(\nu)=\nu_1 \nu_2$. Note that $P$ is such that 
\begin{equation*}
A=P^t \begin{pmatrix} 0&1\\1&0 \end{pmatrix} P.
\end{equation*}
Further note that the choice of $P$ is not unique, since we can multiply by
\begin{equation*}
\pm \begin{pmatrix} \exp(r)&0\\0&\exp(-r) \end{pmatrix}
\end{equation*}
on the left. We take the choice of the sign such that $P^{-1} \left(\begin{smallmatrix} 1\\ -1\end{smallmatrix} \right)\in C_Q$. Now define for $t\in\R$
\begin{equation*}
c(t) := P^{-1} \begin{pmatrix} \exp(t)\\-\exp(-t) \end{pmatrix} \qquad \text{and}\qquad c^\perp(t) := P^{-1} \begin{pmatrix} \exp(t)\\\exp(-t) \end{pmatrix}.
\end{equation*}
We see that $Q(c(t))=-1$, $Q(c^\perp(t))=1$ and $B(c(t),c^\perp (t))=0$. In fact $c$ gives us a parametrization of $C_Q$. If we have $c_i\in C_Q$ then we let $t_i$ be such that
\begin{equation*}
c(t_i) = c_i,
\end{equation*}
and we denote $c^\perp (t_i)$ by $c_i^\perp$. If we make a different choice for $r$, that means that we shift the parameter $t$ by $r$. 

The function that we wish to study is given by
\begin{definition}\label{defphi}
Let $c_1,c_2\in C_Q$ and $a,b\in\R^2$. For $\tau=x+iy\in\Ha$ we define
\begin{equation*}
\begin{split}
\Phi_{a,b} (\tau) &= \Phi_{a,b}^{c_1,c_2} (\tau) \\
&:= \sgn(t_2-t_1)\ y^{1/2} \sum_{\nu\in a+\Z^2} \frac{1}{2}\Bigl[1-\sgn(B(\nu,c_1)B(\nu,c_2))\Bigr] e^{2\pi iQ(\nu)x} e^{2\pi i B(\nu,b)} K_0(2\pi Q(\nu) y)\\
&\ + \sgn(t_2-t_1)\ y^{1/2} \sum_{\nu\in a+\Z^2} \frac{1}{2}\Bigl[1-\sgn(B(\nu,c_1^\perp)B(\nu,c_2^\perp))\Bigr] e^{2\pi iQ(\nu)x} e^{2\pi i B(\nu,b)} K_0(-2\pi Q(\nu) y).
\end{split}
\end{equation*}
\end{definition}
To see that the function is well defined we have to show convergence of the sums, which will be done in the next section.
\begin{remark}
With a little bit of rewriting we see that for $A=\left(\begin{smallmatrix} 3&0\\0&-2\end{smallmatrix}\right)$, $c_1 = \frac{1}{\sqrt{3}}\left(\begin{smallmatrix} -2\\3\end{smallmatrix}\right)$, $c_2 =\frac{1}{\sqrt{3}} \left(\begin{smallmatrix} 2\\3\end{smallmatrix}\right)$, $a=\bigl(\begin{smallmatrix} 1/6\\0\end{smallmatrix}\bigr)$ and $b=\bigl(\begin{smallmatrix} 1/6\\1/4\end{smallmatrix}\bigr)$ we have
\begin{equation*}
\Phi_{a,b}^{c_1,c_2} = \zeta_{12}\ \phi_0,
\end{equation*}
where $\zeta_n$ denotes $e^{2\pi i/n}$ and $\phi_0$ is Cohen's function.
\end{remark}
From the differential equation satisfied by $K_0$ we immediately get
\begin{equation*}
\Delta_0 \Phi_{a,b} = \frac{1}{4} \Phi_{a,b}.
\end{equation*}
However, in general this function does not transform like a modular function. We also consider
\begin{definition}
Let $c_1,c_2\in C_Q$ and $a,b\in\R^2$. For $\tau=x+iy\in\Ha$ we define
\begin{equation*}
\widehat{\Phi}_{a,b} (\tau) = \widehat{\Phi}_{a,b}^{c_1,c_2} (\tau) := y^{1/2} \sum_{\nu\in a+\Z^2} q^{Q(\nu)} e^{2\pi i B(\nu,b)} \int_{t_1}^{t_2} e^{-\pi y B(\nu, c(t))^2} dt,
\end{equation*}
where $q:= e^{2\pi i\tau}$.
\end{definition}
This definition is independent of the choice of $r$. Again convergence of the sum is shown in the next section. From the following theorem we see that this function does transform like a modular function, but in general it is not an eigenfunction of $\Delta_0$ (as we will later see).

\begin{theorem}\label{modtrans}
For $a,b\in\R^2$ we have
\begin{equation*}
\begin{split}
\widehat{\Phi}_{a+\lambda,b+\mu} &= e^{2\pi iB(a,\mu)} \widehat{\Phi}_{a,b} \qquad \text{for all}\ \lambda\in\Z^2\ \text{and}\ \mu\in A^{-1} \Z^2,\\
\widehat{\Phi}_{-a,-b} &= \widehat{\Phi}_{a,b},
\end{split}
\end{equation*}
and the modular transformation properties
\begin{equation*}
\begin{split}
\widehat{\Phi}_{a,b} (\tau+1) &= e^{-2\pi iQ(a)-\pi iB(A^{-1} A^*,a)}\ \widehat{\Phi}_{a,a+b+\frac{1}{2}A^{-1} A^*}(\tau),\\
\widehat{\Phi}_{a,b} (-1/\tau) &= \frac{e^{2\pi i B(a,b)}}{\sqrt{-\det A}} \sum_{p\in A^{-1} \Z^2 \smod{\Z^2}} \widehat{\Phi}_{-b+p,a}(\tau),
\end{split}
\end{equation*}
where $A^*$ is the vector of diagonal elements of $A$.
\end{theorem}
\begin{remark}
From these modular transformation properties we get that if $a,b\in\Q^2$ then $\widehat{\Phi}_{a,b}$ transforms like a modular function on some subgroup of $\operatorname{SL}_2(\Z)$. For convenience we only consider the case that $A$ has integer coefficients. More generally, if we allow the coefficients of $A$ to be rational, we still get that $\widehat{\Phi}_{a,b}$ transforms like a modular function on some subgroup of $\operatorname{SL}_2(\Z)$, but we omit the details.
\end{remark}
So far we have introduced the two functions $\Phi_{a,b}$ and $\widehat{\Phi}_{a,b}$. The first is an eigenfunction of the Laplace operator, but in general not modular. The second is modular, but in general not an eigenfunction of the Laplace operator. We would like to get functions which are both eigenfunctions of the Laplace operator and modular. We will be able to do that in special cases, using the following relation between $\Phi_{a,b}$ and $\widehat{\Phi}_{a,b}$.

\begin{theorem}\label{split}
Let $c_1,c_2\in C_Q$, $b\in\R^2$ and let $a\in\R^2$ be such that $Q$ is non-zero on $a+\Z^2$. Then
\begin{equation*}
\widehat{\Phi}_{a,b}^{c_1,c_2} = \Phi_{a,b}^{c_1,c_2} + \phi_{a,b}^{c_1} - \phi_{a,b}^{c_2},
\end{equation*}
where for $\tau=x+iy\in\Ha$
\begin{equation*}
\phi_{a,b}^{c_0} (\tau) := y^{1/2} \sum_{\nu\in a+\Z^2} \alpha_{t_0}(\nu y^{1/2})\ q^{Q(\nu)} e^{2\pi iB(\nu,b)},
\end{equation*}
with
\begin{equation*}
\alpha_{t_0}(\nu) := \begin{cases} \int_{t_0}^\infty e^{-\pi B(\nu,c(t))^2}dt & \text{if}\ B(\nu,c_0)B(\nu,c_0^\perp)>0,\\
                  -\int_{-\infty}^{t_0} e^{-\pi B(\nu,c(t))^2}dt & \text{if}\ B(\nu,c_0)B(\nu,c_0^\perp)<0,\\
		0 & \text{if}\ B(\nu,c_0)B(\nu,c_0^\perp)=0.
                  \end{cases}
\end{equation*}
\end{theorem}
The main point here is that $\widehat{\Phi}_{a,b}^{c_1,c_2} -\Phi_{a,b}^{c_1,c_2}$ is the difference of two functions, each of which depends only on one $c$. For certain special choices of $a,b$ and $c_1,c_2$ we can have that $\phi_{a,b}^{c_1} = \phi_{a,b}^{c_2}$ and so $\widehat{\Phi}_{a,b}^{c_1,c_2} = \Phi_{a,b}^{c_1,c_2}$ is both an eigenfunction of the Laplace operator and modular. This happens in the case of Cohen's $\phi_0$ function, as we will see in section \ref{example}.

As remarked before, $\widehat{\Phi}_{a,b}^{c_1,c_2}$ is in general not an eigenfunction of the Laplace operator, but instead

\begin{theorem}\label{diff}
We have
\begin{equation*}
\left(\Delta_0 -\frac{1}{4} \right) \widehat{\Phi}_{a,b}^{c_1,c_2} = \frac{\pi}{2} \left(\theta_{a,b}^{c_2} - \theta_{a,b}^{c_1}\right),
\end{equation*}
with
\begin{equation*}
\theta_{a,b}^c (\tau) : = y^{3/2} \sum_{\nu\in a+\Z^2} B(\nu,c) B(\nu,c^\perp) e^{\frac{\pi i}{2} \tau B(\nu,c^\perp)^2 - \frac{\pi i}{2} \overline{\tau} B(\nu,c)^2} e^{2\pi i B(\nu,b)}.
\end{equation*}
\end{theorem}
\begin{remark}
We see that $\theta_{a,b}^c$ ``almost'' splits as the product of a unary theta function of weight 3/2 with the complex conjugate of a unary theta function of weight 3/2. To make this a bit more precise: if $c$ is such that there is an $l\in\R$ for which $lc\in\Z^2$, then from $B(c,c^\perp)=0$ and the assumption that $A$ has integer coefficients we get that there is an $l^\perp \in \R$ for which $l^\perp c^\perp \in\Z^2$. Let $\alpha\in\R$ then be the smallest positive number such that $\alpha \left(\begin{matrix} c& c^\perp \end{matrix}\right)^t A$ has integer coefficients. If we now make the change of variables 
\begin{equation*}
\lambda = \alpha \begin{pmatrix} B(\nu,c)\\ B(\nu,c^\perp)\end{pmatrix} = \alpha \left(\begin{matrix} c& c^\perp \end{matrix}\right)^t A \nu,
\end{equation*}
we get that
\begin{equation*}
\theta_{a,b}^c \in S_{3/2} \otimes y^{3/2} \overline{S_{3/2}}.
\end{equation*}
We leave the details to the interested reader.
\end{remark}

\section{Convergence for $\Phi$ and $\widehat{\Phi}$}\label{conv}
If $c_1=c_2$ (that is $t_1=t_2$) both $\Phi_{a,b}^{c_1,c_2}$ and $\widehat{\Phi}_{a,b}^{c_1,c_2}$ are zero and so we can assume that $t_1\neq t_2$.

We can easily verify that
\begin{equation*}
\begin{split}
c_1 - \exp(t_1-t_2) c_2&= 2 \exp(-t_2) \sinh (t_1-t_2) P^{-1}\begin{pmatrix} 0\\1\end{pmatrix},\\
c_1 - \exp(t_2-t_1) c_2&= 2 \exp(t_2) \sinh (t_1-t_2) P^{-1}\begin{pmatrix} 1\\0\end{pmatrix},
\end{split}
\end{equation*}
and so
\begin{equation}\label{quad}
\begin{split}
&\frac{B(\nu,c_1)^2-2\cosh (t_1-t_2) B(\nu,c_1) B(\nu,c_2) + B(\nu,c_2)^2}{4 \sinh^2(t_1-t_2)} \\
&\qquad \qquad = \frac{\bigl( B(\nu,c_1) - \exp(t_1-t_2) B(\nu,c_2) \bigr) \bigl( B(\nu,c_1) - \exp(t_2-t_1) B(\nu,c_2) \bigr)}{4 \sinh^2(t_1-t_2)}\\
&\qquad \qquad = B\left(\nu , P^{-1}\left(\begin{smallmatrix}0\\1\end{smallmatrix}\right)\right) B\left(\nu , P^{-1}\left(\begin{smallmatrix}1\\0\end{smallmatrix}\right)\right)= Q(\nu),
\end{split}
\end{equation}
where in the last step we have used that $A P^{-1} = P^t \left(\begin{smallmatrix} 0&1\\1&0\end{smallmatrix}\right)$ and $Q(\nu)= (P\nu)_1 (P\nu)_2$.

To get the convergence of the first sum in the definition of $\Phi_{a,b}$, we first observe that $1-\sgn(B(\nu,c_1)B(\nu,c_2))$ is only non-zero if $B(\nu,c_1)B(\nu,c_2)\leq 0$, in which case we get from equation \eqref{quad}
\begin{equation*}
Q(\nu) \geq \frac{B(\nu,c_1)^2 + B(\nu,c_2)^2}{4 \sinh^2(t_1-t_2)},
\end{equation*}
and the right hand side is a positive definite quadratic form. Together with
\begin{equation}\label{est}
0 \leq K_0(x) \leq \sqrt{\frac{\pi}{2x}} e^{-x}
\end{equation}
this shows that the sum converges absolutely.

Completely analogous to equation \eqref{quad} we can show that
\begin{equation*}
-Q(\nu) = \frac{B(\nu,c_1^\perp)^2-2\cosh (t_1-t_2) B(\nu,c_1^\perp) B(\nu,c_2^\perp) + B(\nu,c_2^\perp)^2}{4 \sinh^2(t_1-t_2)},
\end{equation*}
from which we see that on the support of $1-\sgn(B(\nu,c_1^\perp)B(\nu,c_2^\perp))$ we have
\begin{equation*}
-Q(\nu) \geq \frac{B(\nu,c_1^\perp)^2 + B(\nu,c_2^\perp)^2}{4 \sinh^2(t_1-t_2)},
\end{equation*}
which again is positive definite. Using this together with equation \eqref{est} we see that the second sum in the definition of $\Phi_{a,b}$ also
converges absolutely.

From
\begin{equation}\label{Q}
4Q(\nu) = B(\nu,c^\perp(t))^2- B(\nu,c(t))^2,
\end{equation}
which is easy to check, we find
\begin{equation}\label{pos}
Q(\nu) +\frac{1}{2} B(\nu,c(t))^2 = \frac{1}{4} \bigl(B(\nu,c^\perp(t))^2+ B(\nu,c(t))^2\bigr),
\end{equation}
which shows that for all $t\in\R$ the quadratic form $\nu \mapsto Q(\nu) +\frac{1}{2} B(\nu,c(t))^2$ is positive definite, and so
\begin{equation*}
\min_{||\nu||^2 =1} \Bigl(Q(\nu) +\frac{1}{2} B(\nu,c(t))^2 \Bigr)> 0,
\end{equation*}
for all $t\in\R$. If we take $t_1\leq t \leq t_2$ we get that there is an $r\in\R_{>0}$ (which depends on $t_1$ and $t_2$), such that
\begin{alignat*}{2}
Q(\nu) +\frac{1}{2} B(\nu,c(t))^2 & > r, &\qquad\text{for}\ &||\nu||^2 =1,\\
Q(\nu) +\frac{1}{2} B(\nu,c(t))^2 & \geq r ||\nu||^2, &\qquad\text{for}\ &\nu\in\R^2.
\end{alignat*}
So we find that for $t_1 < t_2$
\begin{equation*}
\left| q^{Q(\nu)} \int_{t_1}^{t_2} e^{-\pi y B(\nu,c(t))^2}dt\right| = \int_{t_1}^{t_2} e^{-2\pi y \left(Q(\nu) +\frac{1}{2} B(\nu,c(t))^2\right)} dt \leq (t_2-t_1) e^{-2\pi ry ||\nu||^2},
\end{equation*}
and similarly for $t_1 > t_2$
\begin{equation*}
\left| q^{Q(\nu)} \int_{t_1}^{t_2} e^{-\pi y B(\nu,c(t))^2}dt\right| \leq (t_1-t_2) e^{-2\pi ry ||\nu||^2},
\end{equation*}
which shows the absolute convergence of the sum in the definition of $\widehat{\Phi}_{a,b}$.

\section{Proof of Theorems \ref{modtrans}, \ref{split} and \ref{diff}}
\begin{proof}[Proof of Theorem \ref{modtrans}]
The first two equations are completely trivial and for the modular transformation property with respect to $\tau\mapsto \tau+1$ we observe that $y$ doesn't change and so
\begin{equation*}
\widehat{\Phi}_{a,b} (\tau+1) = y^{1/2} \sum_{\nu\in a+\Z^2} e^{2\pi i Q(\nu)} q^{Q(\nu)} e^{2\pi i B(\nu,b)} \int_{t_1}^{t_2} e^{-\pi y B(\nu, c(t))^2} dt.
\end{equation*}
Using that $A$ has integer coefficients we can easily show that for $\nu\in a+\Z^2$
\begin{equation*}
e^{2\pi iQ(\nu)} = e^{-2\pi i Q(a) -\pi i B(A^{-1} A^* ,a)}\ e^{2\pi iB(\nu,a+\frac{1}{2}A^{-1} A^*)},
\end{equation*}
from which the result follows.

For the modular transformation property with respect to $\tau \mapsto -1/\tau$ we would like to use a theorem by Vign\'eras, found in \cite{vigneras}, which gives a nice general result for indefinite theta functions. However, we can't use the theorem directly, because we need a slightly more general result. The main point is that the function
\begin{equation*}
f_\tau (\nu) = y^{1/2} e^{2\pi iQ(\nu)\tau} \int_{t_1}^{t_2} e^{-\pi yB(\nu,c(t))^2}dt
\end{equation*}
satisfies
\begin{equation}\label{four}
\ef(f_\tau) = \frac{1}{\sqrt{-\det A}}\ f_{-1/\tau}
\end{equation}
(shown below), where $\ef(f)$ is the Fourier transform of $f$, given by
\begin{equation*}
\ef(f)(\nu)= \int_{\R^2} f(\alpha) e^{-2\pi iB(\nu,\alpha)} d\alpha.
\end{equation*}
Note that our normalization of the Fourier transform is different from that used by Vign\'eras. We then use the Poisson summation formula
\begin{equation*}
\sum_{\nu\in\Z^2} f(\nu) = \sum_{\nu\in A^{-1} \Z^2} \ef(f)(\nu)
\end{equation*}
with
\begin{equation*}
f(\nu) = f_\tau (\nu +a) e^{2\pi iB(\nu +a,b)}
\end{equation*}
to get 
\begin{equation*}
\widehat{\Phi}_{a,b} (\tau) = \sum_{\nu\in a+\Z^2} f_\tau(\nu) e^{2\pi i B(\nu,b)} = \sum_{\nu\in\Z^2} f(\nu) = \sum_{\nu\in A^{-1} \Z^2} \ef(f)(\nu).
\end{equation*}
From equation \eqref{four} we then get
\begin{equation*}
\begin{split}
\ef(f)(\nu) &= \int_{\R^2} f_\tau (\alpha +a) e^{2\pi i B(\alpha +a,b)-2\pi iB(\nu,\alpha )}d\alpha\\
&= \ef(f_\tau)(\nu -b) e^{2\pi iB(\nu, a)} = \frac{1}{\sqrt{-\det A}} f_{-1/\tau} (\nu -b)\ e^{2\pi iB(\nu,a)},
\end{split}
\end{equation*}
and so
\begin{equation*}
\widehat{\Phi}_{a,b} (\tau) = \frac{e^{2\pi i B(a,b)}}{\sqrt{-\det A}} \sum_{\nu\in -b +A^{-1} \Z^2} f_{-1/\tau}(\nu)\ e^{2\pi iB(\nu,a)}=\frac{e^{2\pi i B(a,b)}}{\sqrt{-\det A}} \sum_{p\in A^{-1} \Z^2 \smod{\Z^2}} \widehat{\Phi}_{-b+p,a}(-1/\tau).
\end{equation*}
If we replace $\tau$ by $-1/\tau$ we get the desired result. What remains to be shown is that \eqref{four} holds. For this we could use part of the proof of Vign\'eras, by checking that $p(\nu) = \int_{t_1}^{t_2} e^{-\pi B(\nu,c(t))^2}dt$ satisfies a certain differential equation. However, it's fairly easy to get the result directly: we consider
\begin{equation*}
g_\tau (\nu) =  y^{1/2} e^{2\pi iQ(\nu)\tau} e^{-\pi yB(\nu,c(t))^2}= y^{1/2} e^{\frac{\pi i}{2} \tau B(\nu,c^\perp (t))^2 -\frac{\pi i}{2} \overline{\tau} B(\nu,c(t))^2},
\end{equation*}
where the second identity follows from \eqref{Q}. If we make the change of variables
\begin{equation*}
u= \begin{pmatrix} B(\alpha,c^\perp (t))\\ B(\alpha,c(t))\end{pmatrix} = \begin{pmatrix} \exp(-t) & \exp(t)\\ -\exp(-t) & \exp(t)\end{pmatrix} P\alpha,
\end{equation*}
we find
\begin{equation*}
\begin{split}
\ef(g_\tau)(\nu)&= y^{1/2} \int_{\R^2} e^{\frac{\pi i}{2} \tau B(\alpha,c^\perp (t))^2 -\frac{\pi i}{2} \overline{\tau} B(\alpha,c(t))^2-2\pi i B(\nu,\alpha)}d\alpha\\
&= \frac{y^{1/2}}{2|\det P|}\int_\R e^{\frac{\pi i}{2} \tau u_1^2 -\pi i v_1 u_1} du_1 \int_\R e^{-\frac{\pi i}{2} \overline{\tau} u_2^2 -\pi i v_2 u_2} du_2,
\end{split}
\end{equation*}
with
\begin{equation*}
v= \begin{pmatrix} \exp(-t)& \exp(t)\\ \exp(-t) & -\exp(t)\end{pmatrix} P\nu= \begin{pmatrix} B(\nu, c^\perp (t))\\ -B(\nu,c(t)) \end{pmatrix}.
\end{equation*}
Now using $|\det P|=\sqrt{-\det A}$ and
\begin{equation*}
\int_\R e^{\frac{\pi i}{2} \tau u_1^2 -\pi i v_1 u_1} du_1 =\sqrt{\frac{2}{-i\tau}}\ e^{-\frac{\pi i}{2} v_1^2/\tau},\qquad 
\int_\R e^{-\frac{\pi i}{2} \overline{\tau} u_2^2 -\pi i v_2 u_2} du_2=\sqrt{\frac{2}{i\overline{\tau}}}\ e^{\frac{\pi i}{2}v_2^2/\overline{\tau}},
\end{equation*}
we get
\begin{equation*}
\ef(g_\tau) = \frac{1}{\sqrt{-\det A}}\ g_{-1/\tau}.
\end{equation*}
Equation \eqref{four} then follows from $f_\tau = \int_{t_1}^{t_2} g_\tau dt$ by changing the order of integration in the Fourier transform. 
\end{proof}

\begin{proof}[Proof of Theorem \ref{split}]
The proof follows immediately if we use the following lemma with $\nu$ replaced by $\nu y^{1/2}$ and compare the different definitions. Note that each of the sums involved converges absolutely, as seen in section \ref{conv} and the absolute convergence of the sum in the definition of $\phi_{a,b}^{c_0}$ follows from the estimate for $\alpha_{t_0}$ given in the lemma, together with equation \eqref{pos} with $t=t_0$.
\end{proof}

\begin{lemma}\label{lem}
If $Q(\nu) \neq 0$ then
\begin{equation*}
\begin{split}
\int_{t_1}^{t_2} &e^{-\pi B(\nu,c(t))^2}dt = \alpha_{t_1}(\nu) -\alpha_{t_2}(\nu)\ +\\
&\sgn(t_2-t_1) \left( \frac{1}{2} [1-\sgn (B(\nu,c_1) B(\nu, c_2))] +\frac{1}{2} [1-\sgn (B(\nu,c_1^\perp) B(\nu, c_2^\perp))]\right) e^{2\pi Q(\nu)} K_0(2\pi \left| Q(\nu)\right|),
\end{split}
\end{equation*}
and
\begin{equation*}
\left| \alpha_{t_0}(\nu) \right| \leq \frac{e^{-\pi B(\nu,c_0)^2}}{2\sqrt{B(\nu,c_0)^2+B(\nu,c_0^\perp)^2}}.
\end{equation*}
\end{lemma}
\begin{proof}[Proof of Lemma \ref{lem}]
Throughout we assume that $Q(\nu) \neq 0$. First we will show that 
\begin{equation}\label{int}
\int_{-\infty}^{\infty} e^{-\pi B(\nu,c(t))^2}dt = e^{2\pi Q(\nu)} K_0(2\pi |Q(\nu)|),
\end{equation}
and if $B(\nu,c_0)B(\nu,c_0^\perp)=0$ then
\begin{equation}\label{half}
\int_{t_0}^{\infty} e^{-\pi B(\nu,c(t))^2}dt = \frac{1}{2} e^{2\pi Q(\nu)} K_0(2\pi |Q(\nu)|).
\end{equation}
For equation \eqref{int} we observe that
\begin{equation*}
B(\nu,c(t)) = \nu^t A P^{-1} \begin{pmatrix} \exp(t)\\-\exp(-t) \end{pmatrix} = (P\nu)^t \begin{pmatrix} -\exp(-t)\\\exp(t) \end{pmatrix}= -(P\nu)_1 \exp(-t) + (P\nu)_2 \exp(t),
\end{equation*}
and so
\begin{equation*}
B(\nu,c(t))= 2\sgn((P\nu)_2)\sqrt{Q(\nu)} \sinh \left(t+\frac{1}{2} \log \frac{(P\nu)_2}{(P\nu)_1}\right)
\end{equation*}
if $Q(\nu) = (P\nu)_1 (P\nu)_2 > 0$ and
\begin{equation*}
B(\nu,c(t))= 2\sgn((P\nu)_2)\sqrt{-Q(\nu)} \cosh \left(t+\frac{1}{2} \log \frac{-(P\nu)_2}{(P\nu)_1}\right)
\end{equation*}
if $Q(\nu) < 0$, which gives
\begin{equation*}
\int_{-\infty}^{\infty} e^{-\pi B(\nu,c(t))^2}dt = \begin{cases} \int_{-\infty}^\infty e^{-4\pi Q(\nu) \sinh^2 t}dt &\text{if}\ Q(\nu)>0,\\
                                                   \int_{-\infty}^\infty e^{4\pi Q(\nu) \cosh^2 t}dt &\text{if}\ Q(\nu)<0.
                                                   \end{cases}
\end{equation*}
Now using $\sinh^2 t =\frac{1}{2} (-1+\cosh 2t)$ and $\cosh^2 t =\frac{1}{2} (1+\cosh 2t)$ we see that
\begin{equation*}
\int_{-\infty}^{\infty} e^{-\pi B(\nu,c(t))^2}dt = e^{2\pi Q(\nu)} \int_{-\infty}^\infty e^{-2\pi |Q(\nu)| \cosh 2t} dt= e^{2\pi Q(\nu)} K_0(2\pi |Q(\nu)|),
\end{equation*}
where in the last step we used the integral representation
\begin{equation*}
K_0 (x) = \int_{-\infty}^\infty e^{-x \cosh 2t} dt,
\end{equation*}
which follows directly from (see formula 9.6.24 in \cite{stegun})
\begin{equation*}
K_0 (x) = \int_0^\infty e^{-x \cosh t} dt.
\end{equation*}

To prove \eqref{half} we first observe that
\begin{equation*}
c(t_0+t) = c_0 \cosh t + c_0^\perp \sinh t,
\end{equation*}
and so 
\begin{equation}\label{integral}
\int_{t_0}^{\infty} e^{-\pi B(\nu,c(t))^2}dt = \int_{0}^{\infty} e^{-\pi B(\nu,c(t_0+t))^2}dt = \int_{0}^{\infty} e^{-\pi (B(\nu,c_0)\cosh t +B(\nu,c_0^\perp)\sinh t)^2}dt.
\end{equation}
If $B(\nu,c_0^\perp) =0$ this equals
\begin{equation*}
\int_0^\infty e^{-\pi B(\nu,c_0)^2 \cosh^2 t}dt = e^{-\frac{\pi}{2} B(\nu,c_0)^2} \int_0^\infty e^{-\frac{\pi}{2} B(\nu,c_0)^2 \cosh 2t}dt= \frac{1}{2} e^{-\frac{\pi}{2} B(\nu,c_0)^2} K_0\left(\frac{\pi}{2} B(\nu,c_0)^2\right),
\end{equation*}
and if $B(\nu,c_0) =0$ we get
\begin{equation*}
\int_0^\infty e^{-\pi B(\nu,c_0^\perp)^2 \sinh^2 t}dt = e^{\frac{\pi}{2} B(\nu,c_0^\perp)^2} \int_0^\infty e^{-\frac{\pi}{2} B(\nu,c_0^\perp)^2 \cosh 2t}dt= \frac{1}{2} e^{\frac{\pi}{2} B(\nu,c_0^\perp)^2} K_0\left(\frac{\pi}{2} B(\nu,c_0^\perp)^2\right).
\end{equation*}
In both cases we obtain \eqref{half} by using equation \eqref{Q} with $t=t_0$.

Let $\alpha_{t_0}$ be as defined in the Theorem \ref{split}. Using equations \eqref{int} and \eqref{half} we rewrite it as
\begin{equation*}
\alpha_{t_0}(\nu) = \int_{t_0}^\infty e^{-\pi B(\nu,c(t))^2}dt - \frac{1}{2}[1-\sgn(B(\nu,c_0)B(\nu,c_0^\perp))]\ e^{2\pi Q(\nu)} K_0 (2\pi |Q(\nu)|),
\end{equation*}
and so
\begin{equation*}
\begin{split}
\alpha_{t_1}(\nu)-\alpha_{t_2}(\nu) &= \int_{t_1}^{t_2} e^{-\pi B(\nu,c(t))^2}dt\ +\\
& + \frac{1}{2} [\sgn(B(\nu,c_1)B(\nu,c_1^\perp))-\sgn(B(\nu,c_2)B(\nu,c_2^\perp))]\ e^{2\pi Q(\nu)} K_0 (2\pi |Q(\nu)|).
\end{split}
\end{equation*}
If we can show that for $\nu \neq 0$ (note that $Q(\nu)\neq 0$)
\begin{equation}\label{signs}
\begin{split}
-\sgn(B(\nu,c_1)&B(\nu,c_1^\perp))+\sgn(B(\nu,c_2)B(\nu,c_2^\perp))\\
&= \sgn(t_2-t_1) \bigl\{[1-\sgn (B(\nu,c_1) B(\nu, c_2))] + [1-\sgn (B(\nu,c_1^\perp) B(\nu, c_2^\perp))]\bigr\},
\end{split}
\end{equation}
then the result follows. We will prove \eqref{signs} for $t_1<t_2$. The result for $t_1>t_2$ then follows if we interchange $c_1$ and $c_2$. We observe that
\begin{equation*}
\begin{split}
c_2 &= \cosh (t_2-t_1) c_1 + \sinh (t_2-t_1) c_1^\perp,\\
c_2^\perp &= \sinh (t_2-t_1) c_1 + \cosh (t_2-t_1) c_1^\perp,
\end{split}
\end{equation*}
and so
\begin{equation}\label{rel}
\begin{pmatrix} B(\nu,c_2)\\ B(\nu,c_2^\perp)\end{pmatrix} = \begin{pmatrix} \cosh(t_2-t_1) & \sinh(t_2-t_1)\\ \sinh (t_2-t_1) & \cosh (t_2-t_1) \end{pmatrix} \begin{pmatrix} B(\nu,c_1)\\B(\nu,c_1^\perp)\end{pmatrix}.
\end{equation}
If $\nu\neq 0$ then $B(\nu, c_1)$ and $B(\nu,c_1^\perp)$ are not both zero, and since both sides in equation \eqref{signs} are even, it suffices to check it for the case that $B(\nu,c_1),B(\nu,c_1^\perp) \geq 0$ (not both with equality) and the case that $B(\nu,c_1)>0$ and $B(\nu,c_1^\perp)<0$. In the first case we get from equation \eqref{rel} and $t_1<t_2$ that $B(\nu,c_2),B(\nu,c_2^\perp)>0$ and we can verify that equation \eqref{signs} holds. In the second case we get from
\begin{equation*}
\begin{pmatrix} B(\nu,c_1)\\ B(\nu,c_1^\perp)\end{pmatrix} = \begin{pmatrix} \cosh(t_2-t_1) & -\sinh(t_2-t_1)\\ -\sinh (t_2-t_1) & \cosh (t_2-t_1) \end{pmatrix} \begin{pmatrix} B(\nu,c_2)\\B(\nu,c_2^\perp)\end{pmatrix}
\end{equation*}
that we can't have both $B(\nu,c_2)\leq0$ and $B(\nu,c_2^\perp)\geq 0$, since that would mean that $B(\nu,c_1)\leq 0$, and so we have that $B(\nu,c_2)>0$ or $B(\nu,c_2^\perp)<0$. If both cases we can verify that equation \eqref{signs} holds.

To finish the proof we establish the estimate for $\alpha_{t_0}$: assuming $B(\nu,c_0)B(\nu,c_0^\perp)\geq 0$ and $t\geq 0$ we have $\sinh t \geq t$ and
\begin{equation*}
\begin{split}
(B(\nu,c_0)\cosh t &+B(\nu,c_0^\perp)\sinh t)^2\\&= B(\nu,c_0)^2 +\left( B(\nu,c_0)^2+B(\nu,c_0^\perp)^2\right)\sinh^2 t + B(\nu,c_0)B(\nu,c_0^\perp)\sinh 2t\\
&\geq B(\nu,c_0)^2 +\left( B(\nu,c_0)^2+B(\nu,c_0^\perp)^2\right) t^2
\end{split}
\end{equation*}
and so we get from equation \eqref{integral}
\begin{equation*}
\int_{t_0}^{\infty} e^{-\pi B(\nu,c(t))^2}dt \leq e^{-\pi B(\nu,c_0)^2} \int_0^\infty e^{-\pi (B(\nu,c_0)^2+B(\nu,c_0^\perp)^2)t^2} dt= \frac{e^{-\pi B(\nu,c_0)^2}}{2\sqrt{B(\nu,c_0)^2+B(\nu,c_0^\perp)^2}}.
\end{equation*}
Similarly, we find that if $B(\nu,c_0)B(\nu,c_0^\perp)\leq 0$ then
\begin{equation*}
\int_{-\infty}^{t_0} e^{-\pi B(\nu,c(t))^2}dt \leq \frac{e^{-\pi B(\nu,c_0)^2}}{2\sqrt{B(\nu,c_0)^2+B(\nu,c_0^\perp)^2}},
\end{equation*}
which finishes the proof.
\end{proof}

\begin{proof}[Proof of Theorem \ref{diff}]
Using equation \eqref{Q} and
\begin{equation*}
\frac{\partial }{\partial t} c(t) = c^\perp (t), \qquad
\frac{\partial }{\partial t} c^\perp(t) = c(t),
\end{equation*}
it is straightforward to check that 
\begin{equation*}
\left( \Delta_0 -\frac{1}{4} \right) \left\{ y^{1/2} e^{2\pi i Q(\nu)\tau - \pi y B(\nu,c(t))^2} \right\}= \frac{\pi}{2} y^{3/2} \frac{\partial}{\partial t} \left\{B(\nu,c(t)) B(\nu,c^\perp(t)) e^{\frac{\pi i}{2} \tau B(\nu,c^\perp(t))^2 - \frac{\pi i}{2} \overline{\tau} B(\nu,c(t))^2}\right\}.
\end{equation*}
If we use this in the definition of $\widehat{\Phi}_{a,b}^{c_1,c_2}$ we find
\begin{equation*}
\begin{split}
&\left( \Delta_0 -\frac{1}{4} \right) \widehat{\Phi}_{a,b}^{c_1,c_2} (\tau)= \frac{\pi}{2} y^{3/2} \cdot\\
&\cdot \sum_{\nu\in a+\Z^2} \left\{ B(\nu,c_2) B(\nu,c_2^\perp) e^{\frac{\pi i}{2} \tau B(\nu,c_2^\perp)^2 - \frac{\pi i}{2} \overline{\tau} B(\nu,c_2)^2} - B(\nu,c_1) B(\nu,c_1^\perp) e^{\frac{\pi i}{2} \tau B(\nu,c_1^\perp)^2 - \frac{\pi i}{2} \overline{\tau} B(\nu,c_1)^2}\right\} e^{2\pi i B(\nu,b)}\\
&\qquad \qquad = \frac{\pi}{2} \left(\theta_{a,b}^{c_2} (\tau)- \theta_{a,b}^{c_1}(\tau)\right),
\end{split}
\end{equation*}
where in the last step we split the sum into two parts, which is justified, because each part converges individually. This follows from
\begin{equation*}
\left| e^{\frac{\pi i}{2} \tau B(\nu,c^\perp)^2 - \frac{\pi i}{2} \overline{\tau} B(\nu,c)^2}\right|= e^{-\frac{\pi}{2} y (B(\nu,c)^2+B(\nu,c^\perp)^2)}.
\end{equation*}
\end{proof}

\section{Properties of $\phi_{a,b}^c$ and the action of $\operatorname{Aut}^+ (Q,\Z^2)$}

We consider the group of matrices that leave both the quadratic form and the lattice $\Z^2$ invariant, that is
\begin{equation*}
\operatorname{Aut} (Q,\Z^2) =\left\{ \gamma \in \operatorname{GL}_2(\R) \left|\ \gamma^t A \gamma = A,\ \gamma \Z^2=\Z^2\right.\right\}.
\end{equation*}
Let $\gamma$ be an element of this group and let $c\in C_Q$, then $Q(\gamma c) =Q(c)$, so $\gamma C_Q$ is either $C_Q$ or $-C_Q$. Similarly, $\gamma$ can send one of the components of the set of vectors $c\in\R^2$ with $Q(c)=1$ either to itself or the opposite component. We consider only matrices $\gamma$ that leave $C_Q$ invariant, i.e.\ $B(\gamma c,c)<0$ for all $c\in C_Q$, and have determinant 1, so that they also fix the components of $\{c\in\R^2 \mid Q(c)=1\}$. The group of such matrices we denote by $\operatorname{Aut}^+ (Q,\Z^2)$, that is
\begin{equation*}
\operatorname{Aut}^+ (Q,\Z^2):= \left\{ \gamma \in \operatorname{GL}_2(\R) \left|\ Q \circ \gamma = Q,\ \gamma \Z^2=\Z^2,\ \gamma C_Q=C_Q\ \text{and}\ \det(\gamma)=1 \right.\right\}.
\end{equation*}

\begin{lemma}\label{phi}
Let $\phi_{a,b}^c$ be as defined in Theorem \ref{split}. For $c\in C_Q$ and $a,b\in \R^2$ we have
\begin{equation*}
\begin{split}
\phi_{a+\lambda,b+\mu}^c &=e^{2\pi iB(a,\mu)} \phi_{a,b}^c \qquad \text{for all}\ \lambda\in\Z^2\ \text{and}\ \mu\in A^{-1} \Z^2,\\
\phi_{-a,-b}^c &= \phi_{a,b}^c,
\end{split}
\end{equation*}
and
\begin{equation*}
\phi_{\gamma a,\gamma b}^{\gamma c} = \phi_{a,b}^c \qquad \text{for all}\ \gamma\in\operatorname{Aut}^+ (Q,\Z^2).
\end{equation*}
\end{lemma}
\begin{proof}[Proof of the lemma]
The relations given in the first two equations are completely analogous to those given in Theorem \ref{modtrans} and are again trivial. The last part follows directly if we replace $\nu$ by $\gamma\nu$ in the definition of $\phi_{a,b}^c$ and use that
\begin{equation}\label{transalpha}
\alpha_{\gamma t_0} (\gamma \nu) = \alpha_{t_0}(\nu),
\end{equation}
where $\gamma t_0$ denotes the value of the parameter $t$ such that $c(\gamma t_0) = \gamma c(t_0)$. To prove \eqref{transalpha} we note that $c^\perp (\gamma t_0) =\gamma c_0^\perp$ and so $B(\nu,c_0) B(\nu,c_0^\perp)$ doesn't change if we replace both $\nu$ by $\gamma\nu$ and $t_0$ by $\gamma t_0$. Therefore it suffices to show
\begin{equation*}
\begin{split}
\int_{\gamma t_0}^\infty e^{-\pi B(\gamma \nu,c(t))^2} dt &= \int_{t_0}^\infty e^{-\pi B(\nu,c(t))^2} dt,\\
\int_{-\infty}^{\gamma t_0} e^{-\pi B(\gamma \nu,c(t))^2} dt &= \int_{-\infty}^{t_0} e^{-\pi B(\nu,c(t))^2} dt.
\end{split}
\end{equation*}
The first follows directly if we use \eqref{integral} on both sides and the proof of the second is similar.
\end{proof}

\section{An example: $\sigma$ revisited}\label{example}
As an example we take $A=\left(\begin{smallmatrix} 3&0\\0&-2\end{smallmatrix}\right)$ and $c_1 = \frac{1}{\sqrt{3}}\left(\begin{smallmatrix} -2\\3\end{smallmatrix}\right)$, $c_2 =\frac{1}{\sqrt{3}} \left(\begin{smallmatrix} 2\\3\end{smallmatrix}\right)$, and consider
\begin{equation*}
\widehat{\Phi} := \begin{pmatrix} \widehat{\Phi}_{\bigl(\begin{smallmatrix} 1/6\\0\end{smallmatrix}\bigr), \bigl(\begin{smallmatrix} 1/6\\1/4\end{smallmatrix}\bigr)}\\ \sqrt{2}\ \widehat{\Phi}_{\bigl(\begin{smallmatrix} 1/6\\1/4\end{smallmatrix}\bigr), \bigl(\begin{smallmatrix} 1/6\\0\end{smallmatrix}\bigr)} \\ \sqrt{2}\ \widehat{\Phi}_{\bigl(\begin{smallmatrix} 1/6\\1/4\end{smallmatrix}\bigr), \bigl(\begin{smallmatrix} 1/6\\1/4\end{smallmatrix}\bigr)}\end{pmatrix}.
\end{equation*}
Using Theorem \ref{modtrans} we can determine the modular transformation properties of $\widehat{\Phi}$, to find
\begin{equation}\label{mod}
\widehat{\Phi} (\tau +1) = \left(\begin{smallmatrix} \zeta_{24} & 0 & 0\\0&0&\zeta_{48}^5\\0&\zeta_{48}^{-7}&0\end{smallmatrix}\right) \widehat{\Phi}(\tau) \qquad \text{and}\qquad \widehat{\Phi} (-1/\tau) = \left(\begin{smallmatrix} 0&1&0\\1&0&0\\0&0&1\end{smallmatrix}\right) \widehat{\Phi}(\tau),
\end{equation}
where $\zeta_n:= e^{2\pi i/n}$. Note that if we use the theorem for the transformation $\tau \mapsto -1/\tau$, we get the sum of 6 terms (because $A$ has determinant $-6$), but using the relations given in the theorem, together with the extra relation (which we will prove below)
\begin{equation}\label{extra}
\widehat{\Phi}_{\bigl(\begin{smallmatrix}-a_1\\a_2\end{smallmatrix}\bigr),\bigl(\begin{smallmatrix}-b_1\\b_2\end{smallmatrix}\bigr)} = \widehat{\Phi}_{\bigl(\begin{smallmatrix}a_1\\a_2\end{smallmatrix}\bigr),\bigl(\begin{smallmatrix}b_1\\b_2\end{smallmatrix}\bigr)}
\end{equation}
we can see that some of those terms are zero and that the others are multiples of each other.

We can take 
\begin{equation*}
c(t) = \begin{pmatrix} \sqrt{\frac{2}{3}} \sinh t\\ \cosh t\end{pmatrix},
\end{equation*}
which gives $t_1=-t_2$. We then get equation \eqref{extra} by replacing both $\nu_1$ by $-\nu_1$ and $t$ by $-t$ in the definition of $\widehat{\Phi}_{a,b}^{c_1,c_2}$.

From equation \eqref{mod} we see that $\widehat{\Phi}$ transforms like a vector-valued modular function on the full modular group. Using that $\Gamma_0(2)$ is generated by $\left(\begin{smallmatrix} 1&1\\0&1\end{smallmatrix}\right)$ and $\left(\begin{smallmatrix} 1&0\\2&1\end{smallmatrix}\right)=- \left(\begin{smallmatrix} 0&-1\\1&0\end{smallmatrix}\right)\left(\begin{smallmatrix} 1&-2\\0&1\end{smallmatrix}\right)\left(\begin{smallmatrix} 0&-1\\1&0\end{smallmatrix}\right)$ we can easily verify that for the first component of $\widehat{\Phi}$ we have
\begin{equation}\label{mt}
\widehat{\Phi}_{\bigl(\begin{smallmatrix} 1/6\\0\end{smallmatrix}\bigr), \bigl(\begin{smallmatrix} 1/6\\1/4\end{smallmatrix}\bigr)} \left( \frac{a\tau+b}{c\tau+d}\right) = v(\gamma)\ \widehat{\Phi}_{\bigl(\begin{smallmatrix} 1/6\\0\end{smallmatrix}\bigr), \bigl(\begin{smallmatrix} 1/6\\1/4\end{smallmatrix}\bigr)} (\tau)\qquad \text{for all}\ \gamma=\left(\begin{smallmatrix} a&b\\c&d\end{smallmatrix}\right) \in \Gamma_0(2),
\end{equation}
where $v$ is a multiplier system defined uniquely by
\begin{equation*}
v\left(\begin{smallmatrix} 1&1\\0&1\end{smallmatrix}\right) = v\left(\begin{smallmatrix} 1&0\\2&1\end{smallmatrix}\right) = \zeta_{24}.
\end{equation*}
Let $\gamma = \left(\begin{smallmatrix} 5&4\\6&5\end{smallmatrix}\right)$. Then $\gamma \in \operatorname{Aut}^+(Q,\Z^2)$, $c_2=\gamma c_1$ and we can easily verify using Lemma \ref{phi} that
\begin{equation*}
\phi_{\bigl(\begin{smallmatrix} 1/6\\0\end{smallmatrix}\bigr), \bigl(\begin{smallmatrix} 1/6\\1/4\end{smallmatrix}\bigr)}^{c_1} = \phi_{\bigl(\begin{smallmatrix} 1/6\\0\end{smallmatrix}\bigr), \bigl(\begin{smallmatrix} 1/6\\1/4\end{smallmatrix}\bigr)}^{c_2},
\end{equation*}
and so we get from Theorem \ref{split}
\begin{equation*}
\widehat{\Phi}_{\bigl(\begin{smallmatrix} 1/6\\0\end{smallmatrix}\bigr), \bigl(\begin{smallmatrix} 1/6\\1/4\end{smallmatrix}\bigr)} = \Phi_{\bigl(\begin{smallmatrix} 1/6\\0\end{smallmatrix}\bigr), \bigl(\begin{smallmatrix} 1/6\\1/4\end{smallmatrix}\bigr)}=\zeta_{12}\ \phi_0,
\end{equation*}
where $\phi_0$ is Cohen's function. Hence we have reshown that $\phi_0$ is a Maass waveform, with the modular transformation properties explicitly given in \eqref{mt}.

In for example \cite{kane,corson,lovejoy} we find more examples where functions like $\sigma$ and $\sigma^*$ show up. We could consider the corresponding function $\widehat{\Phi}_{a,b}$ and determine its modular transformation behaviour explicitly, like we did here for $\sigma$ and $\sigma^*$. We leave the details to the reader.

\section{A family of Maass waveforms}
We now construct a family of examples were we get Maass waveforms. To get a nice formulation we consider a slightly different version of $\Phi$. For this we consider periodic functions on $\Z^2$, that is, functions $m$ for which there is an $L\in\Z$ such that
\begin{equation*}
m(\nu+\mu) =m(\nu)\qquad \text{for all}\ \nu\in\Z^2\ \text{and}\ \mu\in L\Z^2.
\end{equation*}

\begin{definition}
Let $c_1,c_2\in C_Q$ and let $m$ be a periodic function on $\Z^2$. For $\tau=x+iy\in\Ha$ we define
\begin{equation*}
\begin{split}
\Phi_m^{c_1,c_2} (\tau) &:= \sgn(t_2-t_1)\ y^{1/2} \sum_{\nu\in \Z^2} \frac{1}{2}\Bigl[1-\sgn(B(\nu,c_1)B(\nu,c_2))\Bigr] m(\nu) e^{2\pi iQ(\nu)x} K_0(2\pi Q(\nu) y)\\
&\ + \sgn(t_2-t_1)\ y^{1/2} \sum_{\nu\in\Z^2} \frac{1}{2}\Bigl[1-\sgn(B(\nu,c_1^\perp)B(\nu,c_2^\perp))\Bigr] m(\nu) e^{2\pi iQ(\nu)x} K_0(-2\pi Q(\nu) y).
\end{split}
\end{equation*}
\end{definition}
Then we have
\begin{theorem}
Let $c\in C_Q$ and suppose we have a finite collection $\{ (m_j,\gamma_j)\}$, where $m_j$ is a periodic function on $\Z^2$ such that $Q$ is non-zero on its support, and $\gamma_j\in \operatorname{Aut}^+ (Q,\Z^2)$, such that
\begin{equation*}
\sum_{j} (m_j -m_j\circ \gamma_j) =0,
\end{equation*}
then $\sum_j \Phi_{m_j}^{c,\gamma_j c}$ is a weight 0 Maass form (on some subgroup of $\operatorname{SL}_2(\Z)$, with some multiplier system), with eigenvalue 1/4 of the Laplace operator $\Delta_0$. Further $\sum_j \Phi_{m_j}^{c,\gamma_j c}$ is independent of $c$.
\end{theorem}
\begin{proof}That it's an eigenfunction of the Laplace operator follows from the differential equation for $K_0$. To show that it transforms like a modular function we need
\begin{equation*}
\begin{split}
\widehat{\Phi}_m^{c_1,c_2} (\tau) &:= y^{1/2} \sum_{\nu\in \Z^2} m(\nu)\ q^{Q(\nu)} \int_{t_1}^{t_2} e^{-\pi y B(\nu, c(t))^2} dt,\\
\phi_m^{c_0} (\tau) &:= y^{1/2} \sum_{\nu\in\Z^2} m(\nu)\ \alpha_{t_0}(\nu y^{1/2})\ q^{Q(\nu)}.
\end{split}
\end{equation*}
As in Theorem \ref{split} we have
\begin{equation*}
\widehat{\Phi}_m^{c_1,c_2} = \Phi_m^{c_1,c_2} +\phi_m^{c_1} -\phi_m^{c_2},
\end{equation*}
and as in Lemma \ref{phi}
\begin{equation}\label{inte}
\phi_{m\circ \gamma^{-1}}^{\gamma c} = \phi_m^c,
\end{equation}
and so
\begin{equation*}
\begin{split}
\sum_j \widehat{\Phi}_{m_j}^{c,\gamma_j c} &= \sum_j \Phi_{m_j}^{c,\gamma_j c} + \sum_j \left\{ \phi_{m_j}^c - \phi_{m_j}^{\gamma_j c}\right\}= \sum_j \Phi_{m_j}^{c,\gamma_j c} + \sum_j \left\{ \phi_{m_j}^c - \phi_{m_j \circ \gamma_j}^{c}\right\}\\
&= \sum_j \Phi_{m_j}^{c,\gamma_j c} + \phi_{\sum_j (m_j-m_j\circ \gamma_j)}^c= \sum_j \Phi_{m_j}^{c,\gamma_j c},
\end{split}
\end{equation*}
where we used that $m\mapsto \phi_m^c$ is linear. We can write each $m_j$ as a finite linear combination $\sum_l d_{j,l}\ e^{2\pi iB(\nu,b_{j,l})}$ and so we get $\widehat{\Phi}_{m_j}^{c,\gamma_j c}$ as a finite linear combination $\sum_l d_{j,l} \widehat{\Phi}_{0,b_{j,l}}^{c,\gamma_j c}$. Using Theorem \ref{modtrans} we then see that $\sum_j \widehat{\Phi}_{m_j}^{c,\gamma_j c}=\sum_i \Phi_{m_j}^{c,\gamma_j c}$ transforms like a modular function on some subgroup of $\operatorname{SL}_2(\Z)$. 

Similar to \eqref{inte} we have
\begin{equation*}
\Phi_{m\circ \gamma^{-1}}^{\gamma c_1,\gamma c_2} = \Phi_{m}^{c_1,c_2},
\end{equation*}
for all $\gamma\in\operatorname{Aut}^+(Q,\Z^2)$, and so
\begin{equation*}
\Phi_{m_j}^{c,\gamma_j c} - \Phi_{m_j}^{\overline{c},\gamma_j \overline{c}} = \Phi_{m_j}^{c,\overline{c}} - \Phi_{m_j}^{\gamma_j c,\gamma_j \overline{c}} =\Phi_{m_j}^{c,\overline{c}} - \Phi_{m_j\circ \gamma_j}^{ c,\overline{c}}= \Phi_{m_j-m_j\circ \gamma_j}^{c,\overline{c}},
\end{equation*}
where in the first step we used that $\Phi_m^{c_1,c_2}+\Phi_m^{c_2,c_3} = \Phi_m^{c_1,c_3}$, which follows from \eqref{signs}. If we now sum over all $j$ we get
\begin{equation*}
\sum_j \Phi_{m_j}^{c,\gamma_j c}=\sum_j \Phi_{m_j}^{\overline{c},\gamma_j \overline{c}},
\end{equation*}
which shows the last part of the theorem.
\end{proof}


\begin{thebibliography}{9}
\bibitem{ADH} G. Andrews, F. Dyson and D. Hickerson, \emph{Partitions and indefinite quadratic forms}, Invent.\ Math.\ \textbf{91} (1988), no.\ 3, pages 391--407. 
\bibitem{kane} K. Bringmann and B. Kane, \emph{Multiplicative $q$-hypergeometric series arising from real quadratic fields}, Trans.\ Amer.\ Math.\ Soc., accepted for publication.
\bibitem{cohen} H. Cohen, \emph{$q$-identities for Maass waveforms}, Invent.\ Math.\ \textbf{91} (1988), no.\ 3, pages 409--422.
\bibitem{corson} D. Corson, D. Favero, K. Liesinger and S. Zubairy, \emph{Characters and $q$-series in $\Q(\sqrt{2})$}, J. Number Theory \textbf{107} (2004), pages 392--405.
\bibitem{lovejoy} J. Lovejoy, \emph{Overpartitions and real quadratic fields}, J. Number Theory \textbf{106} (2004), pages 178--186.
\bibitem{stegun} A. Milton and I. Stegun, \emph{Handbook of mathematical functions with formulas, graphs, and mathematical tables}, National Bureau of Standards Applied Mathematics Series, 55 (1964).
\bibitem{vigneras} M.-F. Vign\'eras, \emph{S\'eries th\^eta des formes quadratiques ind\'efinies}, in Modular functions of one variable VI (Proc.\ Second Internat.\ Conf., Univ.\ Bonn, Bonn, 1976), Lecture Notes in Math.\ \textbf{627}, Springer, Berlin (1977), pages 227--239.
\bibitem{zwegers} S. Zwegers, \emph{Mock theta functions}, Ph.D. Thesis, Universiteit Utrecht (2002).
\end{thebibliography}
\end{document}